\input amstex.tex
\documentstyle{amsppt}
\magnification=\magstep1
\vsize 8.5truein
\hsize 6truein

\heading
A Note on the Size of the Largest Ball Inside a Convex Polytope
\endheading
 
\bigskip

\centerline{\bf Imre B\'ar\'any}

\bigskip

\centerline{Alfr\'ed R\'enyi Mathematical Institute,}
\centerline{Hungarian Academy of Sciences}
\centerline{P.O. Box 127, Hungary H-1364}
\centerline{E-mail: barany\@renyi.hu, and}
\centerline{Department of Mathematics}
\centerline{University College London}
\centerline{Gower Street, London, WC1E 6BT, UK} 

\bigskip
 
\centerline{\bf N\'andor Sim\'anyi}

\bigskip

\centerline{University of Alabama at Birmingham}
\centerline{Department of Mathematics}
\centerline{Campbell Hall, Birmingham, AL 35294 U.S.A.}
\centerline{E-mail: simanyi\@math.uab.edu}

\bigskip \bigskip

\hbox{\centerline{\vbox{\hsize 8cm {\bf Abstract.} Let $m>1$ be an integer,
$B_m$ the set of all unit vectors of $\Bbb R^m$ pointing in the direction of
a nonzero integer vector of the cube $[-1,\,1]^m$. Denote by $s_m$ the radius
of the largest ball contained in the convex hull of $B_m$. We determine the
exact value of $s_m$ and obtain the asymptotic equality 
$s_m\sim\frac{2}{\sqrt{\log m}}$.}}}

\bigskip \bigskip

\noindent
Primary subject classification: 52B11

\medskip

\noindent
Secondary subject classification: 52B12

\bigskip \bigskip

\heading
\S1. Introduction
\endheading

\bigskip

Let $m\ge 2$ be an integer, and consider the sets

$$
A_m=\{-1,\,0,\,1\}^m\setminus\{\vec 0\},\text{ and }
B_m=\left\{\frac{v}{||v||}\Bigg|\; v\in A_m\right\}.
$$
Let $C_m$ be the convex hull of $B_m$, and $s_m$ the radius of the
largest ball contained in $C_m$. (Due to the apparent symmetries of $C_m$,
such a largest ball is necessarily centered at the origin.) In the paper
[B-M-S(2005)] (dealing with rotation numbers/vectors of billiards) we needed
sharp lower and upper estimates for the extremal radius $s_m$. Here we
determine the exact value of $s_m$ which, of course, implies such estimates.

\medskip

\subheading{Theorem} 
$$
s_m=\left(\sum_{k=1}^m \frac{1}{\left(\sqrt{k}+\sqrt{k-1}\right)^2}
\right)^{-1/2}
$$

The Theorem implies that 
$$
\frac{1}{4}\log m<s_m^{-2}<\frac{1}{4}\log m+\frac{5}{4}.
$$
As an immediate corollary, the quantity $s_m$ is asymptotically equal to
$\frac{2}{\sqrt{\log m}}$.

\bigskip \bigskip

\heading
\S2. Proof of the Theorem
\endheading

\bigskip

The proof will be split into of a few lemmas. The first one of them is a
trivial observation.

\medskip

\subheading{Lemma 1} The set of vertices $B_m$ of the convex polytope 
$C_m$, and hence $C_m$ itself, is invariant under the action of the full
isometry group $G$ of the cube $[-1, 1]^m$. (The group $G$ is generated by
all permutations of the coordinates in $\Bbb R^m$, and by all reflections
across the coordinate hyperplanes.) \qed

\medskip

We will use the notation $v_k=\frac{1}{\sqrt{k}}\sum_{i=1}^k e_i$ 
($k=1,\dots,m$) for some specific vertices of $C_m$. (Here $e_i$ stands for
the $i$-th standard unit vector of $\Bbb R^m$.)

\medskip

\subheading{Lemma 2} The simplex $S$, spanned by the linearly independent
vectors $v_k$ ($k=1,\dots,m$) as vertices, is a face of the polytope $C_m$
whose outer normal vector is $u=(u_1,\dots,u_m)$ with the coordinates
$u_i=\sqrt{i}-\sqrt{i-1}$.

\medskip

\subheading{Proof} Consider the scalar product function $\langle v,\,u\rangle$
($v\in B_m$) restricted to the set $B_m$ of vertices of the polytope $C_m$.
Elementary inspection shows that this scalar product function can only
attain its maximum value at the vertices $v_k$, and actually,
$$
\langle v_k,\,u\rangle=1
\tag 1
$$
for each $k=1,\dots,m$. This proves all claims of the lemma. \qed

\medskip

\subheading{Lemma 3} For any face $F$ of the polytope $C_m$ there exists
a congruence $g\in G$ such that $g(F)=S$.

\medskip

\subheading{Proof} Fix a non-zero vector $w=(w_1,\dots,w_m)$ whose ray
$R(w)=\left\{\lambda w|\;\lambda\ge 0\right\}$ intersects the interior of
the face $F$. By selecting $w$ in a generic manner, we can assume that
the absolute values $|w_i|$ of its coordinates are distinct and all different 
from zero. Therefore, by applying a suitable element $g\in G$, we can even
assume that 
$$
w_1>w_2>\dots>w_m>0.
\tag 2
$$
We claim that $g(F)=S$. Indeed, by (2) we have the linear expansion

$$
w=\sum_{k=1}^m \sqrt{k}(w_k-w_{k+1})v_k.
$$
of $w$ in the basis $\left\{v_1,\dots,v_m\right\}$ with positive coefficients.
(With the natural convention $w_{m+1}=0$.) This proves that some positive
multiple of $w$ is a convex linear combination of the vertices of $S$ with
non-zero coefficients, so the face $g(F)$ shares an interior point with
$S$. \qed

\medskip

It follows from the previous lemma that the radius $s_m$ of the inscribed
sphere is actually the distance between $S$ and the origin. However, this
distance is equal to $s_m=\langle u,\,e_1\rangle/||u||=1/||u||$ by (1).
It is clear that 
$$
||u||^2=\sum_{k=1}^m \frac{1}{\left(\sqrt{k}+\sqrt{k-1}\right)^2}.
$$
finishing the proof of our theorem. \qed 

Define $R_m=\sum_{k=1}^m\frac{1}{k}$. For the asymptotic value of $s_m$ we use
the elementary fact that $\log m<R_m<\log m+1$.

$$
\frac{1}{4}\log m<\sum_{k=1}^m \frac{1}{4k}<\sum_{k=1}^m
\frac{1}{\left(\sqrt{k}+\sqrt{k-1}\right)^2}=||u||^2
$$

$$
<1+\sum_{k=2}^m\frac{1}{4(k-1)}<1+\frac{1}{4}(\log m+1)=
\frac{1}{4}\log m+\frac{5}{4}.
$$

\medskip

\subheading{Remark 1} Let $K$ be the convex cone generated by the vectors
$v_k$, $k=1,\dots,m$. The meaning of Lemma 3 is that the cones $g(K)$
($g\in G$) form a triangulation of the space $\Bbb R^m$. As a matter of fact,
the intersections of the cones $g(K)$ with the standard $(m-1)$-simplex

$$
S_{m-1}=\left\{x\in\Bbb R^m\big|\;\sum_{i=1}^m x_i=1,\; x_i\ge0\text{ for all }
i\right\}
$$
form the baricentric subdivision of $S_{m-1}$.

\medskip

\subheading{Remark 2} The following natural question has been considered 
in several papers, for instance in [B-F(1988)] and [B-W(2003)]. What is 
the maximal radius $r(m,N)$ of the inscribed ball of the convex hull of
$N$ points chosen from the unit ball of $\Bbb R^m$? In our case $N=3^m-1$ 
and one may wonder how close $s_m$ and $B_m$ are to the maximal radius and 
best arrangement. It turns out that they are very far: it follows from 
the results of [B-F(1988)] and [B-W(2003)] that, 
in the given range $N=3^m-1$, 
$$
r(m,N)=\left( \frac 89 \right)^{1/2}(1+o(1))
$$ 
as $m \to \infty$. So the optimal radius is much larger than $s_m$. This also 
shows that, as expected, $B_m$ is far from being distributed uniformly 
on the unit sphere.

\medskip

\subheading{Acknowledgement} The authors express their sincere gratitude to
Micha\l \newline
Misiurewicz (Indiana University Purdue University at Indianapolis) for
posing the above problem during his joint research with the second author and
Alexander Blokh (University of Alabama at Birmingham.) The first named author
was partially supported by Hungarian National Foundation Grants T 037846 and T
046246.

\bye